\documentclass[12pt]{article}%
\usepackage{graphicx}
\usepackage[intlimits]{amsmath}
\usepackage{latexsym}
\usepackage{amsfonts}
\usepackage{amssymb}%
\makeatletter
\newfont{\footsc}{cmcsc10 at 8truept}
\newfont{\footbf}{cmbx10 at 8truept}
\newfont{\footrm}{cmr10 at 10truept}
\newtheorem{theorem}{Theorem}

\newtheorem{corollary}[theorem]{Corollary}

\newtheorem{proposition}[theorem]{Proposition}

\newenvironment{proof}[1][Proof]{\noindent{\textbf {#1}  }}  {\hfill$\Box$\bigskip}

\begin{document}

\title{Spectral radius and maximum degree of connected graphs}
\author{Vladimir Nikiforov\\Department of Mathematical Sciences, University of Memphis, \\Memphis TN 38152, USA}
\maketitle

\begin{abstract}
Given a connected irregular graph $G$ of order $n,$ write $\mu$ for the
largest eigenvalue of its adjacency matrix, $\Delta$ for its maximum degree,
and $D$ for its diameter. We prove that%
\[
\Delta-\mu>\frac{1}{\left(  D+2\right)  n}%
\]
and this bound is tight up to a constant factor.

This improves previous results of Stevanovi\'{c} and Zhang, and extends a
result of Alon and Sudakov.

\end{abstract}

\section{Introduction}

Our notation is standard (e.g., see \cite{Bol98}, \cite{CDS80}, and
\cite{HoJo88}); in particular, all graphs are defined on the vertex set
$\left\{  1,2,...,n\right\}  =\left[  n\right]  ,$ and $G\left(  n\right)  $
stands for a graph of order $n.$ Given a graph $G=G\left(  n\right)  ,$ we
write $\Delta\left(  G\right)  $ for the maximum degree of $G$ and order the
eigenvalues of the adjacency matrix $A\left(  G\right)  $ as $\mu\left(
G\right)  =\mu_{1}\left(  G\right)  \geq...\geq\mu_{n}\left(  G\right)  $.

Let $G=G\left(  n\right)  $ be an irregular connected graph. It is easy to see
that $\Delta\left(  G\right)  -\mu\left(  G\right)  >0;$ moreover, as
Stevanovi\'{c} \cite{Ste04} proved, in fact
\[
\Delta\left(  G\right)  -\mu\left(  G\right)  \geq\frac{1}{2n\left(
n\Delta-1\right)  \Delta^{2}}.
\]
Zhang \cite{Zha05} improved this, showing that, if $\Delta=\Delta\left(
G\right)  $ and $D$ is the diameter of $G,$ then
\[
\Delta\left(  G\right)  -\mu\left(  G\right)  \geq\frac{1}{\left(
\sqrt{\Delta}+\sqrt{\left(  \Delta-1\right)  }\right)  ^{2}Dn\Delta}\geq
\frac{1}{4Dn\Delta^{2}}\geq\frac{1}{4n\left(  n-1\right)  \Delta^{2}}.
\]
Neither Zhang nor Stevanovi\'{c} give a hint what are the best bounds one can
hope for.

In this note we considerably strengthen these results.

\begin{theorem}
\label{maint}For every connected irregular graph $G$ of order $n$ and diameter
$D,$%
\begin{equation}
\Delta\left(  G\right)  -\mu\left(  G\right)  >\frac{1}{\left(  D+2\right)  n}
\label{mainin}%
\end{equation}
and this bound is tight up to a constant factor.
\end{theorem}

In the same vein, Alon and Sudakov \cite{AlSu01} proved that, if $G$ is a
nonbipartite graph with diameter $D,$ then
\[
\Delta\left(  G\right)  +\mu_{n}\left(  G\right)  >\frac{1}{\left(
D+1\right)  n}%
\]
and this bound is tight up to a constant factor. We extend their result as follows.

\begin{corollary}
For every connected irregular or nonbipartite graph $G$ of order $n$ and
diameter $D,$%
\begin{equation}
\Delta\left(  G\right)  +\mu_{n}\left(  G\right)  >\frac{1}{\left(
D+2\right)  n} \label{mainin1}%
\end{equation}
and this bound is tight up to a constant factor.
\end{corollary}

The rest of the note is organized as follows. In Section \ref{Prel} we present
some preliminary results, including an extension of a theorem of Wei. In
Section \ref{Pf} we give the proof of (\ref{mainin}), and in Section \ref{Ex}
we describe a family of graphs showing the tightness of (\ref{mainin}) and
(\ref{mainin1}).

\section{\label{Prel}Preliminary results}

Given two vectors $\mathbf{x}=\left(  x_{1},...,x_{n}\right)  ,$
$\mathbf{y}=\left(  y_{1},...,y_{n}\right)  ,$ we write $\mathbf{x}%
\leq\mathbf{y}$ to indicate that $x_{i}\leq y_{i}$ for all $i\in\left[
n\right]  .$

Given a graph $G,$ a $k$\emph{-walk} is a sequence of vertices $v_{1}%
,...,v_{k}$ of $G$ such that $v_{i}$ is adjacent to $v_{i+1}$ for all
$i=1,...,k-1;$ we write $w_{k}\left(  G\right)  $ for the number of $k$-walks
in $G$ and $w_{k}\left(  u\right)  $ for the number of $k$-walks starting with
$u.$

Let $G$ be a graph of order $n$ with eigenvalues $\mu_{1}\geq...\geq\mu_{n}$
and $\mathbf{x}_{1},...,\mathbf{x}_{n}$ be unit orthogonal eigenvectors to
$\mu_{1},..,\mu_{n}$. For every $i\in\left[  n\right]  ,$ let $\mathbf{x}%
_{i}=\left(  x_{i1},...,x_{in}\right)  $ and set $\alpha_{i}=\sum_{j=1}%
^{n}x_{ij}.$

The number of $k$-walks in $G$ (see, e.g., \cite{CDS80}, p. 44, Theorem 1.10)
is given as
\begin{equation}
w_{k}\left(  G\right)  =\alpha_{1}^{2}\mu_{1}^{k-1}+...+\alpha_{n}^{2}\mu
_{n}^{k-1} \label{nw}%
\end{equation}

We shall prove the following bound on $\mu\left(  G\right)  $ given in
\cite{Nik05}.

\begin{theorem}
\label{thup}For every integer $k\geq1,$ $r\geq1$ and any graph $G=G\left(
n\right)  ,$
\[
\mu^{r}\left(  G\right)  \leq\max_{u\in V\left(  G\right)  }\frac
{w_{k+r}\left(  u\right)  }{w_{k}\left(  u\right)  }.
\]

\end{theorem}

\begin{proof}
Set $b_{ii}=w_{k}\left(  i\right)  $ for each $i\in\left[  n\right]  $ and let
$B$ be the diagonal matrix with main diagonal $\left(  b_{11},...,b_{nn}%
\right)  .$ Since $B^{-1}A^{r}\left(  G\right)  B$ has the same spectrum as
$A^{r}\left(  G\right)  ,$ the maximum row sum of $B^{-1}A^{r}\left(
G\right)  B$ - say the sum of the $k$th row - is an upper bound on $\mu
^{r}\left(  G\right)  $; hence,
\[
\mu^{r}\left(  G\right)  \leq\sum_{v\in V\left(  G\right)  }w_{r}\left(
k,v\right)  \frac{w_{k}\left(  v\right)  }{w_{k}\left(  k\right)  }%
=\frac{w_{k+r}\left(  k\right)  }{w_{k}\left(  k\right)  }\leq\max_{u\in
V\left(  G\right)  }\frac{w_{k+r}\left(  u\right)  }{w_{k}\left(  u\right)
},
\]
completing the proof.
\end{proof}

Berge \cite{Ber58} attributes the following theorem to Wei \cite{Wei52}.

\begin{theorem}
\label{thWei}If $G$ is a nonbipartite connected graph and $\mathbf{x}=\left(
x_{1},...,x_{n}\right)  $ is the positive unit eigenvector to $\mu\left(
G\right)  ,$ then%
\[
\lim_{k\rightarrow\infty}\frac{1}{w_{k}\left(  G\right)  }\left(  w_{k}\left(
1\right)  ,...,w_{k}\left(  n\right)  \right)  =\frac{1}{x_{1}+...+x_{n}%
}\mathbf{x}.
\]
$\hfill\square$
\end{theorem}

Note that Wei's theorem fails for bipartite connected graphs in general - a
fact overlooked in \cite{Ber58}, \cite{CDS80}, and \cite{CDS95}. Only in
\cite{CRS97}, p. 26, the theorem appears correctly, together with a proof.

Below we extend Wei's theorem for bipartite graphs.

\begin{theorem}
\label{thWei2}Let $G=G\left(  n\right)  $ be a connected bipartite graph and
$\mathbf{x}=\left(  x_{1},...,x_{n}\right)  ,$ $\mathbf{y}=\left(
y_{1},...,y_{n}\right)  $ be unit eigenvectors to $\mu_{1}\left(  G\right)  $
and $\mu_{n}\left(  G\right)  .$ Then%
\[
\lim_{k\rightarrow\infty}\frac{1}{w_{2k+1}\left(  G\right)  }\left(
w_{2k+1}\left(  1\right)  ,...,w_{2k+1}\left(  n\right)  \right)  \geq
\frac{\left\vert x_{1}+...+x_{n}\right\vert -\left\vert y_{1}+...+y_{n}%
\right\vert }{\left\vert x_{1}+...+x_{n}\right\vert ^{2}+\left\vert
y_{1}+...+y_{n}\right\vert ^{2}}\mathbf{x}%
\]
and%
\[
\lim_{k\rightarrow\infty}\frac{1}{w_{2k}\left(  G\right)  }\left(
w_{2k}\left(  1\right)  ,...,w_{2k}\left(  n\right)  \right)  \geq\frac
{1}{\left\vert x_{1}+...+x_{n}\right\vert +\left\vert y_{1}+...+y_{n}%
\right\vert }\mathbf{x}.
\]

\end{theorem}

To prove the theorem we need the following proposition whose proof we omit.

\begin{proposition}
\label{prop1}If $G$ is a bipartite connected graph with parts $A$ and $B,$ and
$\mathbf{x}=\left(  x_{1},...,x_{n}\right)  $ is the positive unit eigenvector
to $\mu_{1}\left(  G\right)  ,$ then $\mu_{n}\left(  G\right)  $ has only two
unit eigenvectors $\mathbf{y}$ and $-\mathbf{y}$ given by
\[
\mathbf{y}=\left(  \varepsilon_{1}x_{1},...,\varepsilon_{n}x_{n}\right)
,\text{ \ where \ }\varepsilon_{i}=\left\{
\begin{array}
[c]{rrr}%
1 & \text{if} & i\in A\\
-1 & \text{if} & i\in B
\end{array}
\right.  \text{.}%
\]
$\hfill\square$
\end{proposition}

\begin{proof}
[\textbf{Proof of Theorem \ref{thWei2} }]Let $G$ be a connected bipartite
graph of order $n$ with eigenvalues $\mu_{1}\geq...\geq\mu_{n}.$ It is known
that
\begin{equation}
\mu_{1}=-\mu_{n},\text{ \ }\left\vert \mu_{i}\right\vert <\mu_{1}\text{ for
}1<i<n. \label{maj}%
\end{equation}
We may and shall assume that $\mathbf{x}$ is a positive vector. Select
orthogonal unit eigenvectors $\mathbf{x}_{1}=\mathbf{x},\mathbf{x}%
_{2},...,\mathbf{x}_{n}=\mathbf{y}$ to $\mu_{1},...,\mu_{n}$.

Letting $\mathbf{i}=\left(  n^{-1/2},...,n^{-1/2}\right)  \in\mathbb{R}^{n}$
and writing $\alpha_{i}$ for the sum of the entries of $\mathbf{x}_{i},$ we
see that $\mathbf{i}=\sum_{i=1}^{n}\alpha_{i}\mathbf{x}_{i}.$ From (\ref{nw})
we see that
\[
\frac{1}{w_{k}\left(  G\right)  }\left(  w_{k}\left(  1\right)  ,...,w_{k}%
\left(  n\right)  \right)  =\frac{A^{k-1}\left(  G\right)  \mathbf{i}}%
{\alpha_{1}^{2}\mu_{1}^{k-1}+...+\alpha_{n}^{2}\mu_{n}^{k-1}}=\frac{\alpha
_{1}\mu_{1}^{k-1}\mathbf{x}_{1}+...+\alpha_{n}\mu_{n}^{k-1}\mathbf{x}_{n}%
}{\alpha_{1}^{2}\mu_{1}^{k-1}+...+\alpha_{n}^{2}\mu_{n}^{k-1}}.
\]
From (\ref{maj}) and Proposition \ref{prop1}, we obtain%
\[
\lim_{k\rightarrow\infty}\frac{\alpha_{1}\mu_{1}^{2k}\mathbf{x}_{1}%
+...+\alpha_{n}\mu_{n}^{2k}\mathbf{x}_{n}}{\alpha_{1}^{2}\mu_{1}%
^{2k}+...+\alpha_{n}^{2}\mu_{n}^{2k}}=\frac{\alpha_{1}\mathbf{x}_{1}%
+\alpha_{n}\mathbf{x}_{n}}{\alpha_{1}^{2}+\alpha_{n}^{2}}\geq\frac{\alpha
_{1}-\left\vert \alpha_{n}\right\vert }{\alpha_{1}^{2}+\alpha_{n}^{2}%
}\mathbf{x},
\]
and%
\[
\lim_{k\rightarrow\infty}\frac{\alpha_{1}\mu_{1}^{2k-1}\mathbf{x}%
_{1}+...+\alpha_{n}\mu_{n}^{2k-1}\mathbf{x}_{n}}{\alpha_{1}^{2}\mu_{1}%
^{2k-1}+...+\alpha_{n}^{2}\mu_{n}^{2k-1}}=\frac{\alpha_{1}\mathbf{x}%
_{1}-\alpha_{n}\mathbf{x}_{n}}{\alpha_{1}^{2}-\alpha_{n}^{2}}\geq\frac
{\alpha_{1}-\left\vert \alpha_{n}\right\vert }{\alpha_{1}^{2}-\alpha_{n}^{2}%
}\mathbf{x}=\frac{1}{\alpha_{1}+\left\vert \alpha_{n}\right\vert }\mathbf{x},
\]
completing the proof.
\end{proof}

\section{\label{Pf}Proof of inequality (\ref{mainin})}

Let $G$ be a graph of order $n$, maximum degree $\Delta,$ largest eigenvalue
$\mu_{1},$ and diameter $D<\infty.$ Let $\mathbf{x}=\left(  x_{1}%
,,...,x_{n}\right)  $ be a positive unit eigenvector to $\mu_{1}.$ Without
loss of generality we shall assume that $x_{1}\geq...\geq x_{n},$ implying, in
particular, that $x_{1}\geq1/\sqrt{n}.$ Set $\alpha_{1}=x_{1}+...+x_{n}$ and
recall that
\[
\alpha_{1}=x_{1}+...+x_{n}\leq\sqrt{n\left(  x_{1}^{2}+...+x_{n}^{2}\right)
}=\sqrt{n}.
\]

Our first goal is to prove that%
\begin{equation}
\mu_{1}\leq\Delta-\frac{1}{\alpha_{1}}x_{n}. \label{crucin}%
\end{equation}

\textbf{Case 1. }$G$ is nonbipartite

Since $\mu_{1}>\left\vert \mu_{i}\right\vert $ for all $1<i\leq n,$ (\ref{nw})
implies that%
\[
\mu_{1}=\lim_{k\rightarrow\infty}\frac{w_{k+1}\left(  G\right)  }{w_{k}\left(
G\right)  }.
\]

Select $j$ to be a vertex with degree $d\left(  j\right)  \leq\Delta-1.$
Clearly, for every $k\geq1,$
\begin{align*}
w_{k+1}\left(  G\right)   &  =\sum_{i\in\left[  n\right]  }w_{k}\left(
i\right)  d\left(  i\right)  =w_{k}\left(  j\right)  d\left(  j\right)
+\sum_{i\in\left[  n\right]  ,i\neq j}w_{k}\left(  i\right)  d\left(
i\right)  \leq-w_{k}\left(  j\right)  +\Delta\sum_{i\in\left[  n\right]
}w_{k}\left(  i\right) \\
&  =\Delta w_{k}\left(  G\right)  -w_{k}\left(  j\right)  .
\end{align*}
Thus, taking the limits, Theorem \ref{thWei} implies that
\[
\mu_{1}=\lim_{k\rightarrow\infty}\frac{w_{k+1}\left(  G\right)  }{w_{k}\left(
G\right)  }\leq\lim_{k\rightarrow\infty}\left(  \Delta-\frac{w_{k}\left(
j\right)  }{w_{k}\left(  G\right)  }\right)  \leq\Delta-\frac{1}{\alpha_{1}%
}x_{j}\leq\Delta-\frac{1}{\alpha_{1}}x_{n},
\]
completing the proof of (\ref{crucin}) for nonbipartite $G.$\bigskip

\textbf{Case 2. }$G$ is bipartite

Select orthogonal unit eigenvectors $\mathbf{x}_{1}>0,\mathbf{x}%
_{2},...,\mathbf{x}_{n}$ to $\mu_{1},...,\mu_{n}$ and write $\alpha_{i}$ for
the sum of the entries of $\mathbf{x}_{i}.$ Since $\mu_{1}=-\mu_{n}$ and
$\mu_{1}>\left\vert \mu_{i}\right\vert $ for all $1<i<n,$ (\ref{nw}) implies
that%
\[
\lim_{k\rightarrow\infty}\frac{w_{2k+1}\left(  G\right)  }{w_{2k}\left(
G\right)  }=\frac{\alpha_{1}^{2}+\alpha_{n}^{2}}{\alpha_{1}^{2}-\alpha_{n}%
^{2}}\mu_{1}.
\]

Select $j$ to be a vertex with degree $d\left(  j\right)  \leq\Delta-1.$
Clearly, for every $k\geq1,$
\[
w_{2k+1}\left(  G\right)  =\sum_{i\in\left[  n\right]  }w_{2k}\left(
i\right)  d\left(  i\right)  \leq-w_{2k}\left(  j\right)  +\Delta\sum
_{i\in\left[  n\right]  }w_{2k}\left(  i\right)  =\Delta w_{2k}\left(
G\right)  -w_{2k}\left(  j\right)  .
\]
Thus, taking the limits, Theorem \ref{thWei2} implies that%
\begin{align*}
\mu_{1}\frac{\alpha_{1}^{2}+\alpha_{n}^{2}}{\alpha_{1}^{2}-\alpha_{n}^{2}}  &
=\lim_{k\rightarrow\infty}\frac{w_{2k+1}\left(  G\right)  }{w_{2k}\left(
G\right)  }\leq\lim_{k\rightarrow\infty}\left(  \Delta-\frac{w_{2k}\left(
j\right)  }{w_{2k}\left(  G\right)  }\right) \\
&  \leq\Delta-\frac{1}{\alpha_{1}+\left\vert \alpha_{n}\right\vert }x_{j}%
\leq\Delta-\frac{1}{\alpha_{1}+\left\vert \alpha_{n}\right\vert }x_{n},
\end{align*}
and so%
\[
\mu_{1}\leq\frac{\alpha_{1}^{2}-\alpha_{n}^{2}}{\alpha_{1}^{2}+\alpha_{n}^{2}%
}\Delta-\frac{\alpha_{1}-\left\vert \alpha_{n}\right\vert }{\alpha_{1}%
^{2}+\alpha_{n}^{2}}x_{n}\leq\Delta-\frac{2\alpha_{n}^{2}}{\alpha_{1}%
^{2}+\alpha_{n}^{2}}\Delta-\frac{\alpha_{1}-\left\vert \alpha_{n}\right\vert
}{\alpha_{1}^{2}+\alpha_{n}^{2}}x_{n}.
\]
Note that
\[
\Delta\geq1\geq\frac{1}{n}\geq\frac{x_{n}}{\alpha_{1}}%
\]
implies
\[
\frac{2\alpha_{n}^{2}}{\alpha_{1}^{2}+\alpha_{n}^{2}}\Delta\geq\left(
\frac{2\alpha_{n}^{2}}{\alpha_{1}\left(  \alpha_{1}^{2}+\alpha_{n}^{2}\right)
}\right)  x_{n}=\frac{1}{\alpha_{1}}x_{n}-\frac{\alpha_{1}-\left\vert
\alpha_{n}\right\vert }{\alpha_{1}^{2}+\alpha_{n}^{2}}x_{n},
\]
completing the proof of (\ref{crucin}).

Having (\ref{crucin}) in hand, the rest of the proof is the same for bipartite
and nonbipartite graphs.

If $x_{n}>\alpha_{1}/\left(  D+2\right)  n,$ the proof is completed, so we
shall assume that
\[
x_{n}\leq\frac{\alpha_{1}}{\left(  D+2\right)  n}\leq\frac{1}{\left(
D+2\right)  \sqrt{n}}.
\]

Since the diameter of $G$ is $D,$ there exists a path $1=v_{1},...,v_{l+1}=n$
of length $l\leq D$ joining $1$ to $n$. Then,%
\begin{align*}
\Delta-\mu_{1}  &  >\sum_{i\in\left[  n\right]  }d\left(  i\right)  x_{i}%
^{2}-\sum_{i\in\left[  n\right]  }\sum_{ij\in E\left(  G\right)  }x_{i}%
x_{j}=\sum_{ij\in E\left(  G\right)  ,i<j}\left(  x_{i}-x_{j}\right)  ^{2}%
\geq\sum_{i=1}^{l}\left(  x_{v_{i}}-x_{v_{i+1}}\right)  ^{2}\\
&  \geq\frac{1}{l}\left(  \sum_{i=1}^{l}\left\vert x_{v_{i}}-x_{v_{i+1}%
}\right\vert \right)  ^{2}\geq\frac{1}{l}\left(  \sum_{i=1}^{l}x_{v_{i}%
}-x_{v_{i+1}}\right)  ^{2}\geq\frac{1}{D}\left(  x_{1}-x_{n}\right)  ^{2}\\
&  >\frac{1}{D}\left(  \frac{1}{\sqrt{n}}-\frac{1}{\left(  D+2\right)
\sqrt{n}}\right)  ^{2}>\frac{1}{Dn}\left(  1-\frac{2}{\left(  D+2\right)
}\right)  =\frac{1}{\left(  D+2\right)  n},
\end{align*}
completing the proof.$\hfill\square$

\section{\label{Ex}A family of bipartite graphs}

In this section, for all $\Delta\geq2$ and $k\geq1,$ we shall contsruct an
irregular bipartite graph $G_{\Delta,k}$ of order $n=2k\Delta,$ maximum degree
$\Delta\left(  G_{\Delta,k}\right)  =\Delta,$ and diameter $D=4k-1$ such that
\[
\Delta\left(  G_{\Delta,k}\right)  -\mu_{1}\left(  G_{\Delta,k}\right)
<\frac{\pi^{2}}{2\Delta\left(  k+1\right)  ^{2}}<\frac{4\pi^{2}}{nD}.
\]
Therefore, the graphs $G_{\Delta,k}$ show that inequalities (\ref{mainin}) and
(\ref{mainin1}) are tight up to constant factor.

Select $k$ disjoint copies of the complete bipartite graph $K_{\Delta,\Delta
},$ say $G_{1},...,G_{k}.$ Remove an edge, say $v_{2i-1}v_{2i},$ from each
$G_{i}.$ Join $v_{2i}$ to $v_{2i+1}$ for each $i=1,...,k-1$ and write
$G_{\Delta,k}$ for the resulting graph. Clearly, the order of $G_{\Delta,k}$
is $2k\Delta,$ $\Delta\left(  G_{\Delta,k}\right)  =\Delta,$ and its diameter
is $4k-1.$

\begin{proposition}%
\[
\mu_{1}\left(  G_{\Delta,k}\right)  >\Delta-\frac{\pi^{2}}{2\Delta\left(
k+1\right)  ^{2}}.
\]

\end{proposition}

\begin{proof}
We shall use the Rayleigh principle. Assuming $V\left(  G_{\Delta,k}\right)
=\left[  2\Delta k\right]  ,$ for all $j=1,...,k$ and $i\in V\left(
G_{j}\right)  ,$ set
\[
x_{i}=\frac{1}{\sqrt{\left(  k+1\right)  \Delta}}\sin\frac{j\pi}{k+1}.
\]

Recall that the largest eigenvalue $\mu_{1}\left(  P_{k}\right)  $ of the path
$P_{k}$ of order $k$ is $\mu_{1}\left(  P_{k}\right)  =2\cos\frac{\pi}{n+1}$
and the unit vector
\[
\sqrt{\frac{2}{k+1}}\left(  \sin\frac{\pi}{k+1},\sin\frac{2\pi}{k+1}%
,...,\sin\frac{k\pi}{k+1}\right)
\]
is an eigenvector to $\mu_{1}\left(  P_{k}\right)  $ (see, e.g., \cite{CRS97},
p. 27.) Hence,%
\begin{align*}
\mu_{1}\left(  G_{\Delta,k}\right)   &  \geq\sum_{ij\in E\left(  G\right)
}x_{i}x_{j}=\sum_{l=1}^{k}\sum_{ij\in E\left(  G_{l}\right)  }x_{i}x_{j}%
+2\sum_{i=1}^{k}x_{v_{2i}}x_{v_{2i+1}}\\
&  =\frac{\Delta^{2}-1}{\Delta}\sum_{l=1}^{k}\frac{2}{k+1}\sin^{2}\frac{2l\pi
}{k+1}+\frac{1}{\Delta}\cos\frac{\pi}{n+1}\\
&  =\frac{\Delta^{2}-1}{\Delta}+\frac{1}{\Delta}\cos\frac{\pi}{n+1}%
>\Delta-\frac{1}{\Delta}+\frac{1}{\Delta}\cos\frac{\pi}{k+1}\\
&  =\Delta-\frac{\pi^{2}}{2\Delta\left(  k+1\right)  ^{2}},
\end{align*}
completing the proof.
\end{proof}

.

\end{document}